\documentstyle[12pt,fullpage]{article}
\begin{document}

\newlength{\originalbase}
\setlength{\originalbase}{\baselineskip}
\newcommand{\spacing}[1]{\setlength{\baselineskip}{#1\originalbase}}

\newenvironment{proof}{\noindent{\bf\boldmath Proof}:\,\,} {\hfill
$\Box$\medskip} 

\def\f{{\mathbf f\,}}
\def\g{{\mathbf g\,}}
\def\F{{\mathcal F}}
\def\0{{\mathbf 0}}
\def\1{{\mathbf 1}}
\def\y{{\mathbf y\,}}
\def\x{{\mathbf x\,}}
\def\deg{{\mathtt deg}}
\def\diag{{\mathtt diag}}
\def\bsl{{\!\,-\,\!}}
\def\T{{\mathsf ^T}}
\def\rk{{\mathrm rk}\,}
\def\sign{{\mathrm sign}\,}
\def\d{{o\!\ell}}
\def\R{\hbox{\sf I\kern-.12em R}}
\def\sR{\hbox{\scriptsize\sf I\kern-.12em R}}
\def\N{\hbox{\sf I\kern-.12em N}}
\def\Q{\hbox{\sf C\kern-.47em Q}}
\def\C{\hbox{\sf C\kern-.47em C}}
\def\sC{\hbox{\scriptsize\sf C\kern-.47em C}}
\def\proofend{\hfill$\Box$

\medskip}

\begin{center}
{\LARGE\sf
The Coppersmith-Tetali-Winkler Identity for Mechanical Systems
}\\[8mm]
{\sc Andre\u{\i} Kotlov}\\[3.5mm]
{\small CWI, Kruislaan 413, 1098 SJ Amsterdam, the Netherlands}\\
{\tt andrei}@{\tt cwi.nl}\\[8mm]
\end{center}

\begin{abstract}
{\small
\baselineskip=0.158in
Using a mixture of linear algebra and statics, we derive what can be
viewed as a slight generalization of the Coppersmith-Tetali-Winkler
Identity $H(i,j)+H(j,k)+H(k,i)=H(j,i)+H(k,j)+H(i,k)$ for hitting times
of a random walk.  
}
\end{abstract}
\section{Linear Algebra}

In this paper, the bold-faced characters are reserved for denoting
vectors; the dimension of the vectors is $n\ge 3$, except when the context
suggests otherwise.  For example, $\f$ is a vector with co\"ordinates
$f_1,\dots,f_n$, while $\0$ and $\1$ are the all-zero and all-one
vectors, respectively.  

Consider the system of linear equations 
$$
(1)\qquad A\x=\f
$$
where $\x$ is the vector of the unknowns, $\f$ is a fixed vector, and $A$
is 
an $n\times n$ matrix such that:
\begin{enumerate}
\item[(i)]{$A$ is symmetric;}
\item[(ii)]{the off-diagonal entries of $A$ are non-positive;}
\item[(iii)]{$A$ is irreducible, which is to say that the simple 
graph whose
adjacency matrix has the same off-diagonal zero pattern as $A$ is
connected (throughout the paper, we refer to this graph as ``{\em the 
underlying graph $G$}'');}
\item[(iv)]{$A\1=\0$.}
\end{enumerate}
Applied to a matrix with properties (i)--(iii), the
Perron-Frobenius Theorem states that its least eigenvalue  
has multiplicity one, and the corresponding eigenvector can be
chosen to have all its co\"ordinates positive.  Hence, by (iv), $A$ 
is positive semidefinite and of rank $n-1$.  Consequently, (1)
has no solution if $\f\cdot\1\not=0$ and a solution unique up to a
shift by a multiple of $\1$ if $\f\cdot\1=0$.  For
$i\in\{1,\dots,n\}$, set $\f_{i}:=(f_1,\dots,f_i-\f\cdot\1,\dots, f_n)$
(so that $\f_i\cdot\1=0$) and let $\x_i=(x_{i1},\dots,x_{in})$ denote
the [unique] solution to the system of linear equations
$$
(2)\qquad A\x=\f_i
$$
such that $x_{ii}=0$.  Our goal is to demonstrate the identity
$$
(3) \qquad x_{12}+x_{23}+x_{31}=x_{21}+x_{32}+x_{13}.
$$
Suppose first that $n>3$ and let us apply [symmetric] Gaussian
elimination to (2), to eliminate the non-zero off-diagonal entries of
$A$ in the last column and row.  
The reader can easily see that the first $n-1$ unknowns are not
affected by such an elimination, while the upper left
$(n-1)\times(n-1)$ corner, $A'$, of the matrix obtained from $A$, and
the $(n-1)$-dimensional vector, $\f_i'$, whose entries are the same
as the first $n-1$ entries of the vector obtained from $\f_i$, have
the same properties as $A$ and $\f_i$, respectively.  
In other words, the system of $n-1$ linear equations 
$$
(4)\qquad A'\x'=\f_i'
$$
with $n-1$ unknowns $x_1,\dots,x_{n-1}$ can be considered in lieu of
(2), and if $\x_i=(x_{i1},\dots,x_{in})$ is a solution to (2)
then $\x_i':=(x_{i1},\dots,x_{i,n-1})$ is a solution to (4).  Hence,
in order to establish (3), we may assume that $n=3$.

Still, to verify (3) for $n=3$ directly is a rather tedious task.
After all, one would have to solve {\em three} systems of three
equations each {\em in general}.
Instead, we prefer the following detour. In the next section, we will
see that (3) is trivially true if the underlying graph $G$ is the star 
$K_{1,n-1}$.  If $n=3$ then $G$ is either the star $K_{1,2}$ or the
triangle $K_3$.  In the latter case, the matrix $A$ has the form
$$
\left[\begin{array}{ccc}\alpha_2+\alpha_3&-\alpha_3&-\alpha_2\\
-\alpha_3&\alpha_1+\alpha_3&-\alpha_1\\
-\alpha_2&-\alpha_1&\alpha_1+\alpha_2\end{array}\right]
$$
where 
$\alpha_1,\alpha_2,$ and $\alpha_3$ are some positive numbers.  Set
$c:=\alpha_1\alpha_2+\alpha_2\alpha_3+\alpha_3\alpha_1$ and 
$\beta_i:=c/\alpha_i$, $i\in\{1,2,3\}$.  Then, the matrix
$$
B:=\left[\begin{array}{cccc}    \beta_1&0&0&-\beta_1\\
				0&\beta_2&0&-\beta_2\\
				0&0&\beta_3&-\beta_3\\
-\beta_1&-\beta_2&-\beta_3&\beta_1+\beta_2+\beta_3
\end{array}\right]
$$
has properties (i)--(iv) and it is not difficult to verify that $A$ is
the upper left $3\times 3$ corner of the matrix obtained from $B$ by
symmetric Gaussian elimination.  Hence, (3) holds for $A$ if it holds
for $B$.  Which it does, since the underlying graph for
$B$ is the star $K_{1,3}$.

\section{Mechanical Model}
The systems (1) and (2) have the following mechanical interpretation:
$n$ masses labelled $1$ through $n$ are 
submerged in water, so that each mass $j$ is subject to a (positive or 
negative or zero) weight, $f_j$.  In addition, every two masses
$j$ and $k$ are connected by a rubber band of resistance
$a_{jk}=a_{kj}\ge 0$ where $-a_{jk}$ is the $jk$ entry of $A$.  (A
rubber band is of resistance at most $a\ge 0$ if and only if
a force of magnitude $a\ell$ is sufficient to stretch it to length
$\ell$.)    
If mass $i$ is nailed to the origin, such a system of masses has a
unique equilibrium.  Let $x_{ij}$ denote the 
altitude of mass $j$ in this equilibrium. (So that $x_{ii}=0$.)
The fact that the system is in an equilibrium is expressed by 
equating the sum of the forces acting on every mass $j$ to zero.
In an equilibrium, these forces are the weight $f_j$, the reaction 
$-\f\cdot\1=-(f_1+\dots+f_n)$ of the nail if $j=i$, and the 
resistance forces $a_{1j}(x_{i1}-x_{ij}),\dots,a_{nj}(x_{in}-x_{ij})$.
In other words, the vector $(x_{i1},\dots,x_{in})$ is the unique
solution $\x_i$ to the system (2) such that $x_{ii}=0$.

Using this physical interpretation of (2), we see that (3) is
trivially true if the underlying graph $G$ is
the star $K_{1,n-1}$ with say center $k$.  This is because
$x_{ij}=x_{ik}+x_{kj}$ in this case.  This completes the proof of (3).

\medskip

\noindent{\bf Remark 1:\,}
In the mechanical model, we could replace the weights $f_i$ by
arbitrary $d$-dimensional forces and, correspondingly, the altitudes
$x_{ij}$ by the $d$-dimensional position vectors.  Then, the
``$d$-dimensional generalization'' of (3) would be proved
co\"ordinate-wise.

\medskip

\noindent{\bf Remark 2:\,}
The identity $x_{i_1i_2}+x_{i_2i_3}+\dots+x_{i_ki_1}
=x_{i_2i_1}+x_{i_3i_2}+\dots+x_{i_1i_k}$ can be shown for any sequence
$i_1,\dots,i_k$ of indices using (3); we leave this to the reader.

\section{Random Walks}

Let $G$ be a connected graph on $n\ge 3$ vertices.  
To obtain a {\bf random walk} out of vertex $j$ in $G$, one
recursively builds a sequence $i_0i_1i_2\dots$ of vertices, in
which $i_0=j$ and, for $t\ge 0$, the vertex $i_{t+1}$ is chosen among
the neighbors of $i_t$ with the uniform probability.  
The {\bf hitting time $H(j,i)$} is the expectation of the least number
$t$ such that $i_t=i$. (So that $H(i,i)=0$.)  If $i\not=j$ then, clearly,
$$
(5)\qquad H(j,i)=1+{1\over\deg(j)}\,\sum_{\{k\,:\,jk\in G\}}H(k,i).
$$
Multiplying both sides of (5) by $\deg(j)$,
we see that  
the vector $\x_i:=(H(1,i),\dots,H(n,i))$ is a solution to the system
of linear equations
$$
(6)\qquad A\x=\g_i
$$ 
where $\g_i$ is {\em some} vector
coinciding with the vector  $\f:=(\deg(1)\dots,\deg(n))$ in all the
co\"ordinates except perhaps $i$, and $A$ is the difference between
the diagonal matrix $\diag\,\{\f\}$ and 
the adjacency matrix of $G$.  Clearly,  $A$ satisfies
(i)--(iv), whence $\g_i$ {\em must\,} be the vector $\f_i$
obtained from $\f$ as in Section~1.  Hence (3), taking in this context
the form of the 
Coppersmith-Tetali-Winkler Identity~\cite{1} 
$H(i,j)+H(j,k)+H(k,i)=H(j,i)+H(k,j)+H(i,k)$, holds.  Also,
substituting $\x_i$ into equation $i$ of (6) gives the identity
$$
\sum_{\{k\,:\,ik\in G\}}H(k,i)=2m-\deg(i)
$$
where $2m:=\deg(1)+\dots+\deg(n)$ is twice the number of edges in $G$.
It follows that the {\bf return time} $R(i)$, i.e.~the expectation of
the smallest {\em natural\,} $t$ such that $i_t=i$ in a random walk
out of $i$, is
$$
1+{1\over\deg(i)}(2m-\deg(i))={2m\over\deg(i)}.
$$

\bigskip

\noindent{\Large \bf Acknowledgements}

\bigskip

\noindent

I thank L\'aszl\'o Lov\'asz, Tam\'as Fleiner, and Peter Winkler
for their help. 

\end{document}